\documentclass[a4paper]{amsart}

\usepackage[utf8]{inputenc}
\usepackage[english]{babel}

\usepackage{amsthm}
\usepackage{amsmath}
\usepackage{amssymb}

\usepackage{color}
\usepackage{hyperref}
\usepackage{enumitem}
\setlist{itemsep=4pt plus 2pt minus 1pt,font=\normalfont}

\usepackage{tikz}
\usepackage{url}

\newcommand{\down}[1]{{\downarrow}\,#1} 
\newcommand{\up}[1]{{\uparrow}\,#1}
\newcommand{\outn}[1]{#1^{+}}
\newcommand{\inn}[1]{#1^{-}}
\newcommand{\oa}{\outn{\alpha}}
\newcommand{\ob}{\outn{\beta}}
\newcommand{\ia}{\inn{\alpha}}
\newcommand{\ib}{\inn{\beta}}
\newcommand{\po}[1]{T(#1)}
\newcommand{\poa}{\po{\alpha}}
\newcommand{\park}[1]{Q_{#1}}

\newcommand{\mat}[2]{\mathbf{M}_{#1}(#2)}
\newcommand{\zero}{\mathbf{0}}
\newcommand{\one}{\mathbf{1}}
\newcommand{\two}{\mathbf{2}}
\newcommand{\boole}[1]{\mat{#1}{\mathbf{2}}}

\newcommand{\set}{\Omega}
\newcommand{\elem}{\omega}
\newcommand{\cutt}[1]{\Gamma_{#1}}
\newcommand{\cut}{\cutt{\elem}}
\newcommand{\cutk}{\cutt{k}}
\newcommand{\pow}{\mathcal{P}}
\newcommand{\x}{\mathbf{x}}

\theoremstyle{plain}
\newtheorem{theorem}{Theorem}[section]
\newtheorem{lemma}[theorem]{Lemma}
\newtheorem{proposition}[theorem]{Proposition}
\newtheorem{corollary}[theorem]{Corollary}

\theoremstyle{definition}
\newtheorem{definition}[theorem]{Definition}
\newtheorem{example}[theorem]{Example}

\newtheorem{remark}[theorem]{Remark}

\begin{document}

\title{Multiplication of matrices over lattices}

\author[K. K\'{a}tai-Urb\'{a}n]{Kamilla K\'{a}tai-Urb\'{a}n}
\address[K. K\'{a}tai-Urb\'{a}n]{Bolyai Institute, 
University of Szeged, 
Aradi v\'{e}rtan\'{u}k tere 1, 
H-6720 Szeged, Hungary}
\email[K. K\'{a}tai-Urb\'{a}n]{katai@math.u-szeged.hu}

\author[T.~Waldhauser]{Tam\'{a}s Waldhauser}
\address[T. Waldhauser]{Bolyai Institute, 
University of Szeged, 
Aradi v\'{e}rtan\'{u}k tere 1, 
H-6720 Szeged, Hungary}
\email[T.~Waldhauser]{twaldha@math.u-szeged.hu}


\date{\today}

\begin{abstract}
We study the multiplication operation of square matrices over lattices. 
If the underlying lattice is distributive, then matrices form a semigroup; we investigate idempotent and nilpotent elements and the maximal subgroups of this matrix semigroup.
We prove that matrix multiplication over nondistributive lattices is antiassociative, and we determine the invertible matrices in the case when the least or the greatest element of the lattice is irreducible.
\end{abstract}


\maketitle

\section{Introduction}\label{sec:intro}

Matrix multiplication of matrices over a lattice $L$ can be defined in the same way as for matrices over rings, letting the join operation play the role of addition and the meet operation play the role of multiplication. 
For notational convenience, we will actually write the lattice operations as addition and multiplication.
Thus, throughout the paper, $L=(L;+,\cdot)$ denotes an arbitrary lattice, and $\mat{n}{L}$ stands for the set of all $n\times n$ matrices over $L$.
To exclude trivial cases, we will always assume without further mention that $L$ has at least two elements and $n\geq2$.
If $L$ has a least and a greatest element (these will be denoted by $0$ and $1$), then we can define the identity matrix $I\in\mat{n}{L}$ with ones on the diagonal and zeros everywhere off the diagonal, and it is easy to see that $I$ is indeed the identity element of $\mat{n}{L}$.

Matrix multiplication is not always associative, and if it is not, then we may ask how far it is from being associative.
There are several ways to measure associativity; one of them is the associative spectrum introduced in \cite{CsaWal2000}.
The number of possibilities of inserting parentheses (or brackets) into a product $x_1\cdot\ldots\cdot x_n$ is given by the $(n-1)$-st Catalan number $C_{n-1}=\frac{1}{n-1}\binom{2n-2}{n-1}$.
If multiplication is associative, then all these different \emph{bracketings} give the same result, but if the multiplication is not associative, then some of the bracketings may induce different $n$-variable term functions.
The \emph{associative spectrum} of a binary operation is the sequence $\{s_n\}_{n=1}^{\infty}$ that counts the number of different term functions induced by bracketings of the product $x_1\cdot\ldots\cdot x_n$.
Clearly $s_1=s_2=1$, and $1\leq s_n\leq C_{n-1}$ holds for all natural numbers $n$, and we can say that the faster the spectrum grows, the less associative the multiplication is. 
In particular, if the associative spectrum is the sequence of Catalan numbers, then the multiplication is said to be \emph{antiassociative}.
Of course, there are plenty of operations that fall between the two extreme cases of being associative or antiassociative; examples of associative spectra of various growth rates can be found in \cite{CsaWal2000,LieWal2009}.

We shall see in Section~\ref{sec:arbi} that there is a dichotomy for matrix multiplication over lattices: if $L$ is distributive, then $\mat{n}{L}$ is a semigroup, while if $L$ is not distributive, then the multiplication of $\mat{n}{L}$ is antiassociative.
Nonassociativity has some unfortunate consequences: powers of matrices, nilpotent matrices and inverse matrices are not always well-defined.
On the other hand, we prove that if $L$ is bounded and at least one of $0$ and $1$ is irreducible, then inverses are unique (even if $L$ is not distributive), and we describe explicitly the invertible matrices in this case, showing that they form a group isomorphic to the symmetric group $S_n$.

In Section~\ref{sec:boole} we focus on the semigroup $\boole{n}$ of $n\times n$ matrices over the two-element lattice $\two=\{0,1\}$.
We can regard a matrix $A\in\boole{n}$ as the characteristic function of a set $\alpha\subseteq X^2$ where $X:=\{1,\dots,n\}$, thus matrices over $\two$ correspond to binary relations, and $\boole{n}$ is isomorphic to the semigroup of binary relations on the set $X$.
We present various results about this semigroup in Section~\ref{sec:boole}, namely, we describe its idempotent elements, Green's relations $\mathcal{D}$ and $\mathcal{H}$ for idempotents, and the maximal subgroups ``around" the idempotents.
Some of these results appear already in the literature, but we provide elementary proofs for the readers' convenience, offering a gentle introduction to the structure of the semigroup $\boole{n}$ for non-semigroup theorists.

Section~\ref{sec:distri} deals with matrices over bounded distributive lattices. 
Boundedness is not a serious restriction, since most of the time we shall work in a finitely generated sublattice (for instance, in the sublattice generated by the $n^{2}$ entries of an $n\times n$ matrix), and finitely generated distributive lattices are finite. 
A matrix $A\in\mat{n}{L}$ can be viewed as a multiple-valued analogue of a binary relation $\alpha\subseteq X^2$.
Generalizing results of Section~\ref{sec:boole} to this multiple-valued setting, we describe idempotents and maximal subgroups around some special idempotents in $\mat{n}{L}$; the full description of maximal subgroups constitutes a topic for further research.
We also determine nilpotent matrices over distributive lattices with a meet-irreducible bottom element (this includes chains, which are the most important cases from the viewpoint of applications), and then apply it to solve a problem arising from applications of fuzzy relations \cite{SesTep2011}.

Some personal remarks from the second author about Ivo Rosenberg: 
As a graduate student working in clone theory under the supervision of B\'{e}la Cs\'{a}k\'{a}ny, I certainly learned the name of Ivo Rosenberg early in my studies. 
His theorems on maximal and minimal clones are cornerstones of the theory clones, and I always imagined the discoverer of these theorems as an unapproachable ``giant". 
It is no wonder that I was thrilled to meet him at the AAA58 conference in Vienna in 1999.
Unfortunately, it was our first and last personal encounter.
We spoke only a few words, and he apologized very kindly for not being able to attend my talk. 
I was a bit disappointed, but much more astonished, for receiving such friendly apologies from this giant of clone theory as a first-year doctoral student.
My talk was about measuring associativity, and our joint paper with B\'{e}la Cs\'{a}k\'{a}ny about associative spectra appeared in this journal 20 years ago, in the special issue dedicated to the 65th birthday of Ivo Rosenberg.
Now this is a special issue for a much more sad occasion, and I can only hope that this modest contribution is worthy to commemorate Ivo Rosenberg.

\section{Matrices over arbitrary lattices}\label{sec:arbi}
\subsection{Antiassociativity of matrix multiplication}
First we characterize lattices with associative matrix multiplication.
Here, and in the rest of the paper, we always assume that all matrices
are square matrices of size at least $2 \times 2$.

\begin{proposition}
\label{prop mult. ass. <==> L dist.}Multiplication of matrices over 
a lattice $L$ is associative if and only if $L$ is a distributive lattice.
\end{proposition}

\begin{figure}
\centering
\begin{tabular}{lcr}
\begin{tikzpicture}[scale=1];
\draw[line width=0.2mm] (1,0) node [below]  {$0$} -- (2,1)
node [right]  {$c$};
\draw[line width=0.2mm] (1,0)  -- (0,1) node [left] {$a$};
\draw[line width=0.2mm] (1,0)  -- (1,1) node [left] {$b$};
\draw[line width=0.2mm] (1,1)  -- (1,2)
node [above]  {$1$};
\draw[line width=0.2mm] (0,1) -- (1,2);
\draw[line width=0.2mm] (2,1) -- (1,2);
\draw[fill=white] (1,0) circle (.1cm);
\draw[fill=white] (0,1) circle (.1cm);
\draw[fill=white] (2,1) circle (.1cm);
\draw[fill=white] (1,1) circle (.1cm);
\draw[fill=white] (1,2) circle (.1cm);
\end{tikzpicture}
& 
&
\begin{tikzpicture}[scale=1];
\draw[line width=0.2mm] (1,0) node [below]  {$0$} -- (2,0.5)
node [right]  {$a$};
\draw[line width=0.2mm] (1,0)  -- (0,1) node [left] {$b$};
\draw[line width=0.2mm] (2,0.5)  -- (2,1.5) node [right] {$c$};
\draw[line width=0.2mm] (0,1)  -- (1,2)
node [above]  {$1$};
\draw[line width=0.2mm] (2,1.5) -- (1,2);
\draw[fill=white] (1,0) circle (.1cm);
\draw[fill=white] (2,0.5) circle (.1cm);
\draw[fill=white] (2,1.5) circle (.1cm);
\draw[fill=white] (0,1) circle (.1cm);
\draw[fill=white] (1,2) circle (.1cm);
\end{tikzpicture}
\end{tabular}

\caption{The lattices $M_3$ and $N_5$}
\label{fig:M3N5}
\end{figure}
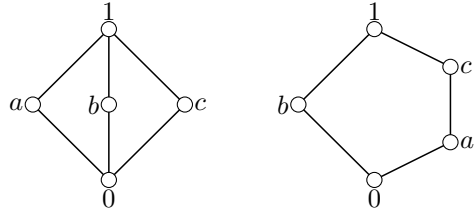

\begin{proof}
If $L$ is distributive, then one can prove associativity of matrix
multiplication in the same way as it is proved for matrices over commutative
rings. In fact, all the usual properties of matrix operations hold in this
case (e.g., multiplication is distributive over addition).

If $L$ is not distributive, then $M_{3}$ or $N_{5}$ embeds into $L$ (see Figure~\ref{fig:M3N5}), so it
suffices to prove nonassociativity of matrix multiplication over these two
lattices.  
Let us consider the following three matrices from $\mat{2}{M_3}$ or from $\mat{2}{N_5}$:%
\[
A=%
\begin{pmatrix}
a & b\\
0 & 0
\end{pmatrix}
,\quad B=%
\begin{pmatrix}
1 & 0\\
1 & 0
\end{pmatrix}
,\quad C=%
\begin{pmatrix}
c & 0\\
0 & 0
\end{pmatrix}
.
\]
Then it is easy to verify that $\left(  AB\right)  C\neq A\left(  BC\right)$.
For any $n \geq 2$, we can construct matrices attesting the
nonassociativity of multiplication in $\mat{n}{L}$ by inserting $A$, $B$ and $C$
into the top left $2 \times 2$ corner of an $n \times n$ matrix and filling 
all the remaining entries with $0$.
\end{proof}

We can strengthen Proposition~\ref{prop mult. ass. <==> L dist.}; if $L$ is 
not distributive, then matrix multiplication over $L$ is not merely nonassociative:
it is antiassociative! We derive this as a corollary of the following general
proposition.
\begin{proposition}\label{prop:nonass and unit ==> antiass} 
If  a binary operation has an identity element, then it is either associative
(i.e., the associative spectrum is constant $1$) or it is antiassociative
(i.e., the associative spectrum consists of the Catalan numbers).
\end{proposition}
\begin{proof}
Let $(G;\cdot)$ be a groupoid with an identity element $1$, and assume that
$(ab)c \neq a(bc)$ for some $a,b,c \in G$. We prove by induction on $n$ that
any two bracketings $p \neq q$ of size $n$ induce different term operations 
on $G$. For $n=1,2$ this claim is void, and for $n=3$ it holds by the 
nonassociativity of the multiplication of $G$.
Assume now that different bracketings of size less than $n$ induce different
term functions, and let $p,q$ be two distinct bracketings of size $n$.

First we consider the case when the ``outermost" multiplication of $p$ and $q$ 
is at the same place: $p=p_1(x_1,\dots,x_k) \cdot p_2(x_{k+1},\dots,x_n)$
and $q=q_1(x_1,\dots,x_k) \cdot q_2(x_{k+1},\dots,x_n)$.
Since $p$ and $q$ are not the same term, we have $p_1 \neq q_1$ or $p_2 \neq q_2$
(perhaps both). If $p_1 \neq q_1$, then, by the induction hypothesis, 
there exist elements $a_1,\dots,a_k \in G$ such that 
$p_1(a_1,\dots,a_k) \neq q_1(a_1,\dots,a_k)$. This implies
\begin{align*}
p(a_1,\dots,a_k,1,\dots,1) &= p_1(a_1,\dots,a_k) \cdot p_2(1,\dots,1) \\
&= p_1(a_1,\dots,a_k) \cdot 1 = p_1(a_1,\dots,a_k) \\
&\neq q_1(a_1,\dots,a_k) = q_1(a_1,\dots,a_k) \cdot 1  \\
&= q_1(a_1,\dots,a_k) \cdot q_2(1,\dots,1) \\
&= q(a_1,\dots,a_k,1,\dots,1),   
\end{align*}
thus the term functions corresponding to $p$ and $q$ are indeed different.
If $p_2 \neq q_2$, then a similar argument can be used, assigning the value
$1$ to the variables $x_1,\dots,x_k$.

Now assume that the outermost multiplications in $p$ and $q$ are not
at the same place: $p=p_1(x_1,\dots,x_k) \cdot p_2(x_{k+1},\dots,x_n)$
and $q=q_1(x_1,\dots,x_\ell) \cdot q_2(x_{\ell+1},\dots,x_n)$, where
$k \neq \ell$. We may suppose without loss of generality that $k<\ell$.
Let us put $x_1=a$, $x_{k+1}=b$, $x_{\ell+1}=c$, and assign the value $1$
to all the remaining variables. Then $p$ evaluates to
\[
p_1(a,1,\dots,1) \cdot p_2(b,1,\dots,1,c,1,\dots,1) = a(bc),
\]
while $q$ gives the value
\[
q_1(a,1,\dots,1,b,1,\dots,1)\cdot q_2(c,1,\dots,1) = (ab)c,
\]
proving that $p$ and $q$ induce different term functions, as claimed.
\end{proof}

\begin{corollary}
If the lattice $L$ is not distributive, then the multiplication
of matrices over $L$ is antiassociative.
\end{corollary}
\begin{proof}
Since $L$ is not distributive, it has a sublattice $K$
that is isomorphic to $M_3$ or to $N_5$. The lattice $K$ is 
bounded, hence $\mat{n}{K}$ has an identity element, 
thus its multiplication is antiassociative by 
propositions~\ref{prop mult. ass. <==> L dist.} and
\ref{prop:nonass and unit ==> antiass}.
This implies antiassociativity of the multiplication of
$\mat{n}{L}$, as it contains $\mat{n}{K}$ as a subgroupoid.
\end{proof}

The following example shows that for nondistributive
lattices even the definition of a power of a matrix 
and the notion of nilpotence can be problematic.

\begin{example}
Let $A$ be the following $5 \times 5$ matrix over $M_3$:
\[
A=
\begin{pmatrix}
0 & a & 0 & 0 & 0\\
0 & 0 & b & c & 0\\
0 & 0 & 0 & 0 & b\\
0 & 0 & 0 & 0 & c\\
0 & 0 & 0 & 0 & 0
\end{pmatrix}.
\]
Then we have $(AA)A = 0 \neq A(AA)$. Thus $A$ has two different
``cubes"; one of them is zero, the other one is not.
\end{example}

\subsection{Invertible matrices}
As another illustration of the unpleasant consequences of
nonassociativity, we present an example of a matrix 
having several inverses.

\begin{example}
Consider the following two matrices over $N_5$:
\[
A=%
\begin{pmatrix}
c & b\\
b & c
\end{pmatrix}
,\quad B=%
\begin{pmatrix}
a & b\\
b & c
\end{pmatrix}
.
\]
Then we have $AA=AB=BA=I$, thus $A$ and $B$ are both inverses of $A$.
\end{example}

Let us conclude this section with some positive results:
we will prove that if $L$ is a bounded lattice such
that at least one of $0$ and $1$ is irreducible, then
inverses in $\mat{n}{L}$ are unique; moreover, the
only invertible matrices are the permutation matrices.
For any permutation $\pi \in S_n$, we define the 
\emph{permutation matrix} corresponding to $\pi$
as the matrix 
$P_\pi=(p_{ij})_{i,j=1}^{n}\in\mat{n}{L}$
given by
\[
p_{ij}=
\begin{cases}
1, & \text{if $j=\pi(i)$;}\\
0, & \text{otherwise.}
\end{cases}
\]
\begin{remark}\label{rem:permutation matrices}
Just as over commutative rings, the matrix $P_\pi A$
is obtained from $A$ by permuting its rows according
to the permutation $\pi$; similarly, $A P_\pi$
is obtained from $A$ by permuting its columns according
to the permutation $\pi^{-1}$. In particular, 
we have $P_\pi P_\sigma = P_{\pi\sigma}$ for all
$\pi,\sigma \in S_n$, and the (unique) inverse
of $P_\pi$ is $P_{\pi^{-1}}$.
\end{remark}

\begin{theorem}\label{thm:01irr invertible}
Let $L$ be a bounded lattice in which $0$ (the bottom element) is meet-irreducible or $1$ (the top element) is join-irreducible.
Then for all matrices $A,B \in \mat{n}{L}$, we have $AB=I$ if and only if $A=P_\pi$ and $B=P_{\pi^{-1}}$ for some permutation $\pi \in S_n$.
\end{theorem}

\begin{proof}
The ``if" part is clear 
(see Remark~\ref{rem:permutation matrices}); so
we only prove the ``only if" part. 
Moreover, it suffices to prove that $A=P_\pi$; then
$B=P_{\pi^{-1}}$ follows by 
Remark~\ref{rem:permutation matrices}.
First we make some general observations, assuming only
boundedness about the lattice $L$.

Let $A,B \in \mat{n}{L}$ with $AB=I$.
Considering the diagonal entries of $AB=I$,
we have $\sum_{j=1}^n a_{ij}b_{ji}=1$ for all 
$i=1,\dots,n$. This implies that for each $i$ there is 
at least one $j$ such that $a_{ij}b_{ji} \neq 0$. 
Denoting such an index  $j$ by $\pi(i)$, we get a map
$\pi\colon \{ 1,\dots,n \} \to \{ 1,\dots,n\}$
such that
\begin{equation}\label{eq:diagonal}
a_{i\pi(i)} \neq 0 \text{ and } b_{\pi(i)i} \neq 0
\text{ for all $i \in \{ 1,\dots,n \}$}.
\end{equation}
The off-diagonal entries of $AB=I$ yield
$\sum_{j=1}^n a_{ij}b_{jk}=0$ whenever $i \neq k$,
hence
\begin{equation}\label{eq:off-diagonal}
a_{ij}b_{jk}=0
\text{ for all $i,j,k \in \{ 1,\dots,n \}$ 
with $i \neq k$}.
\end{equation}

Assume first that $1$ is join-irreducible. 
Then at least one of the summands in 
$\sum_{j=1}^n a_{ij}b_{ji}=1$ must be $1$, 
hence we can replace \eqref{eq:diagonal} by
the following stronger condition:
\begin{equation}\label{eq:diagonal1}
\tag{\ref{eq:diagonal}'}
a_{i\pi(i)}=b_{\pi(i)i}=1
\text{ for all $i \in \{ 1,\dots,n \}$}.
\end{equation}
Now we can see that $\pi$ is injective:
if we had $\pi(i)=\pi(k)=:j$ for some 
$i \neq k$, then \eqref{eq:diagonal1} would imply that 
$a_{ij}=b_{jk}=1$, contradicting
\eqref{eq:off-diagonal}.
In order to prove that $A$ is a permutation matrix,
let us consider an entry $a_{ij}$ in $A$ with
$j \neq \pi(i)$. Letting $k=\pi^{-1}(j)$,
we have $b_{jk}=1$ by \eqref{eq:diagonal1};
on the other hand, \eqref{eq:off-diagonal} implies $a_{ij}b_{jk}=0$, as $i \neq k$. Thus $a_{ij}=0$
whenever $j \neq \pi(i)$, and this together with 
\eqref{eq:diagonal1} proves that $A=P_\pi$. 

Suppose next that $0$ is meet-irreducible.
Then \eqref{eq:off-diagonal} takes the following form:
\begin{equation}\label{eq:off-diagonal0}
\tag{\ref{eq:off-diagonal}'}
a_{ij}=0 \text{ or } b_{jk}=0
\text{ for all $i,j,k \in \{ 1,\dots,n \}$ 
with $i \neq k$}.
\end{equation}
Again, $\pi$ is injective:
if we had $\pi(i)=\pi(k)=:j$ for some 
$i \neq k$, then \eqref{eq:diagonal} would imply that 
$a_{ij} \neq 0$ and $b_{jk} \neq 0$, contradicting
\eqref{eq:off-diagonal0}.
Just as in the previous case, we can prove that
$a_{ij}=0$ whenever $j \neq \pi(i)$. Indeed,
for $k=\pi^{-1}(j)$
we have $b_{jk} \neq 0$ by \eqref{eq:diagonal},
and then \eqref{eq:off-diagonal0} implies $a_{ij}=0$.
To show that $A=P_\pi$, it only remains to prove
that $a_{i\pi(i)}=1$ for every $i$. This follows from
the following inequality:
\[
1 = \sum_{j=1}^n a_{ij}b_{ji} = a_{i\pi(i)}b_{\pi(i)i}
\leq a_{i\pi(i)}.
\]
\end{proof}

\begin{remark}
As a consequence of
Theorem~\ref{thm:01irr invertible},
we have that $AB=I$ implies $BA=I$ for all matrices 
$A,B \in \mat{n}{L}$ if $L$ satisfies the irreducibility 
condition of the theorem.
For semigroups (and also for rings), the property
$AB=I \implies BA=I$ is called \emph{Dedekind-finiteness}.
\end{remark}

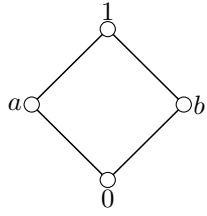
\begin{figure}
\centering
\begin{tikzpicture}[scale=1];
\draw[line width=0.2mm] (1,0) node [below]  {$0$} -- (2,1)
node [right]  {$b$};
\draw[line width=0.2mm] (1,0)  -- (0,1) node [left] {$a$};
\draw[line width=0.2mm] (2,1)  -- (1,2)
node [above]  {$1$};
\draw[line width=0.2mm] (0,1) -- (1,2);
\draw[fill=white] (1,0) circle (.1cm);
\draw[fill=white] (0,1) circle (.1cm);
\draw[fill=white] (2,1) circle (.1cm);
\draw[fill=white] (1,2) circle (.1cm);
\end{tikzpicture}
\caption{The lattice $\two\times\two$.}
\label{fig:2x2}
\end{figure}

\begin{example}
Theorem~\ref{thm:01irr invertible} is not necessarily valid if neither $0$ nor $1$ is irreducible. 
As an example, let 
$A=\begin{pmatrix}
a&b\\
b&a
\end{pmatrix}$
over the lattice $\two\times\two$ shown in Figure~\ref{fig:2x2}.
This lattice is distributive, hence $\mat{n}{L}$ is a semigroup and inverses are unique.
It is easy to verify that $A$ has an inverse (in fact, we have $A^{-1}=A$), even though $A$ is not a permutation matrix.
\end{example}

\section{Matrices over the two-element chain}\label{sec:boole}

To each $n \times n$ matrix $A$ over the two-element lattice
$\two = \{ 0,1 \}$, we can associate a binary relation $\alpha$ 
defined on the set $X := \{ 1,\dots,n \}$, 
by letting $(i,j) \in \alpha \iff a_{ij}=1$.
Matrix multiplication translates to relational product
in this interpretation: if the relations
corresponding  to $A,B \in \boole{n}$ are $\alpha$ and $\beta$, 
then $AB$ describes the relation $\alpha \circ \beta$.
Thus $\boole{n}$ is isomorphic to the semigroup of binary 
relations on the $n$-element set. This semigroup plays a prominent 
role in semigroup theory; we recall a few of the plethora
of results about this semigroup in this section, and we
also recast some of the proofs in a simple form.

\begin{remark}\label{rem:trucks}
We can also regard the relation $\alpha \subseteq X^2$ 
corresponding  to $A \in \boole{n}$ as the edge set of a 
directed graph with vertex set $X$, having $A$ as its adjacency 
matrix. We can think of this graph as a transportation network:
the vertices are sites, and the edges are (possibly one-way)
roads, on which trucks can transport goods between the sites.
If $a_{ii}=0$ (i.e., $(i,i) \notin \alpha$), then trucks
are not allowed to stop at site  $i$, while if
$a_{ii}=1$ (i.e., there is a loop $(i,i) \in \alpha$),
then there is a parking lot at site $i$, where trucks can
wait as long as they wish.
Powers of $A$ account for 
routes\footnote{We use the term \emph{route} 
for a sequence of connecting edges (with 
possible repetitions). The usual terminology
would be \emph{walk}, but we would like
to avoid the uncanny image of a walking truck...} 
in our graph: if 
$A^\ell=(w_{ij})_{i,j=1}^{n}$, then
$w_{ij}=1$ if and only if there is a directed 
route of length $\ell$ from $i$ to $j$.
\end{remark}

\subsection{Idempotent matrices}
The characterization of  idempotent elements of $\boole{n}$ was 
given by B.~Schein \cite{Schein1970} in terms of so-called 
pseudo-orders.
A reflexive transitive relation is called a \emph{quasi-order}.
The symmetric part $\alpha \cap \alpha^{-1}$ of a quasi-order
$\alpha$ is an equivalence relation, and $\alpha$ induces
a natural partial order on the blocks of this equivalence relation.
We usually use the symbol $\leq$ for a quasi-order on the set $X$,
and we denote the corresponding equivalence relation by $\sim$.
Thus the partially ordered set (poset, for short) corresponding
to the quasi-order $\leq$ is $(X/\!\sim;\leq)$.
We say that an element $y \in X$ \emph{covers} $x \in X$
(notation: $x \prec y$), if $x/\!\sim$ is strictly less than 
$y/\!\sim$, and there is no third $\sim$-block between them:
\[
x \prec y \iff x \leq y, x \nsim y \text{ and }
\forall z \in X\colon x \leq z \leq y \implies 
x \sim z \text{ or } z \sim y.
\]

A pseudo-order relation is obtained from a quasi-order by
removing some of the loops (i.e., edges of the form $(x,x)$)
in such a way, that loops can be removed only from singleton 
$\sim$-blocks, and it is not allowed to remove loops from
both members of a covering pair.
\begin{definition}\label{def:pseudo-order}
Let $\alpha \subseteq X^2$ be a binary relation, and let
$\park{\alpha}$ denote the set of vertices with a loop:
$\park{\alpha}=\{x\in X : (x,x)\in\alpha\}$.
We say that $\alpha$ is a \emph{pseudo-order} if the reflexive closure 
$\alpha \cup \{ (x,x) : x \in X\setminus \park{\alpha}\}$
is a quasi-order (we denote this quasi-order by $\leq$
and we use the symbols $\sim$ and $\prec$ for the 
corresponding equivalence relation and cover relation),
and $\park{\alpha}$ satisfies the following two conditions:
\begin{enumerate}[label=(\alph*)]
\item $\forall x \in X\setminus \park{\alpha}\colon\ 
x/\!\sim\, = \{ x \}$,
\item $\forall x,y \in X\colon\  x \prec y \implies 
x \in \park{\alpha} \text{ or } y \in \park{\alpha}$.
\end{enumerate}
\end{definition}

\begin{remark}\label{rem:pseudo truck}
Let us give an interpretation of pseudo-orders in terms
of the transportation network outlined in Remark~\ref{rem:trucks}.
A relation $\alpha \subseteq X^2$ is a pseudo-order
if and only if whenever you drive from site $x$ to site $y$,
\begin{enumerate}[label=(\alph*')]
\item you can choose a direct route (formally: if there
is a route from $x$ to $y$, then $(x,y)$ is an edge), and
\item it is also possible to plan your route so that you will 
have a chance to take a rest in a parking lot on the way
(formally: if there is a route from $x$ to $y$, then
there is a route that includes a vertex with a loop).
\end{enumerate}
Indeed, condition (a) in Definition~\ref{def:pseudo-order} 
ensures that the removal of loops from the underlying 
quasi-order $\leq$ does not ruin transitivity, thus (a') holds
for every pseudo-order. (Observe that (a') is actually
equivalent to transitivity.) 
To verify (b'), choose a longest possible
route that does not pass through any $\sim$-block more than once;
then each edge in this route is a covering pair, and
at least one member of a covering pair has a loop 
(if any of them belongs to a non-singleton 
$\sim$-block, then  condition (a), otherwise condition (b)
provides a loop).

Conversely, let us assume that (a') and (b') hold for $\alpha$, 
and let us denote the reflexive closure of $\alpha$ by $\leq$.
Condition (a') implies that $\alpha$ is transitive, hence $\leq$ 
is also transitive, thus it is a quasi-order. Transitivity of $\alpha$ also implies that (a) holds. To verify (b), consider a
covering pair $x \prec y$. 
If the $\sim$-block of $x$ or $y$ is not a singleton, then
condition (a) shows that there is a loop at $x$ or $y$.
Otherwise, by the definition of covering,
no route from $x$ to $y$ passes through any vertex other
than $x$ and $y$. Therefore, the parking lot guaranteed
by (b') must be at $x$ or at $y$, and this proves (b).
\end{remark}

Now we are ready to state and prove the characterization
of idempotent binary relations given by Schein \cite{Schein1970}.
We will use the description of pseudo-orders given
in Remark~\ref{rem:pseudo truck}.

\begin{theorem}\label{thm:idemp iff pseudo}
A matrix over $\two$ is idempotent if and only if 
the corresponding binary relation is a pseudo-order.
\end{theorem} 
\begin{proof}
Let $\alpha$ be the binary relation on $X$ corresponding to
the matrix $A \in \boole{n}$. As a preliminary observation,
let us note that $\alpha$ is transitive if and only if
$\alpha \circ \alpha \subseteq \alpha$, which 
in turn is equivalent to $A^2 \leq A$.

Assume first that $A$ is idempotent. Then $A^2 \leq A$,
so $\alpha$ is transitive, hence  condition (a') 
of Remark~\ref{rem:pseudo truck} holds. Idempotence
of $A$ implies $A=A^2=A^3=\dots$, thus whenever there
is a route from $x$ to $y$, there are arbitrarily long routes 
from $x$ to $y$.
A long enough route must include a directed cycle, and
every vertex of such a cycle has a loop, by transitivity.
This proves (b'), therefore $\alpha$ is a pseudo-order.

Now let us suppose that $\alpha$ is a pseudo-order.
Then $\alpha$ is transitive by condition (a'), hence
$A^2 \leq A$. Multiplying this inequality by $A^{m-1}$, we get
$A^{m+1} \leq A^m$ for every positive integer $m$, thus
the powers of $A$ form a decreasing sequence:
$A \geq A^2 \geq A^3 \geq \cdots$. Since $\boole{n}$ is a finite
set, this sequence cannot be strictly decreasing, i.e., there
is a positive integer $\ell$ such that
\begin{equation}\label{eq:convergent sequence}
A \geq A^2 \geq A^3 \geq \dots 
\geq A^{\ell} = A^{\ell+1} = A^{\ell+2} = \dots 
= \lim_{m\rightarrow\infty}A^{m}.
\end{equation}
Here the limit is understood in the discrete topology on
$\boole{n}$, but this is not very important, as an ultimately 
constant sequence converges in every topology.
For every edge $(x,y) \in \alpha$, condition (b') provides
a route from $x$ to $y$ with a parking lot on the way. 
We can park there as long as we wish, before continuing 
our trip to $y$, thus there are arbitrarily long routes 
from $x$ to $y$.
This means that $A \leq \lim_{m\rightarrow\infty}A^{m}$, 
which together with the inequalities of 
\eqref{eq:convergent sequence} implies that
$A = A^2 = A^3 = \dots$, hence $A$ is idempotent.
\end{proof}

\subsection{Green's relations}
If $e$ is an idempotent element in a semigroup, then
there is a maximal subgroup $H_e$ ``around" $e$, having
$e$ as its identity element. Having determined the idempotents
of $\boole{n}$, our next goal is to describe 
the maximal subgroups corresponding to these idempotents.
We will need Green's equivalence relations $\mathcal{L}$,
$\mathcal{R}$, $\mathcal{H}$ and $\mathcal{D}$, which can 
be defined in any semigroup, but we write out the definition
for the semigroup $\mat{n}{L}$, where $L$ is a distributive 
lattice. Two elements $A,B\in \mat{n}{L}$ are in $\mathcal L$ 
relation if  they generate the same principal left ideal, that 
is, if and only if  there exist $C,D\in\mat{n}L$ such that 
$CA=B$ and $DB=A$. Similarly, the relation $\mathcal R$ can be 
defined by  $(A,B)\in \mathcal R $ if and only if there exist 
$C,D\in\mat{n}L$ such that $AC=B$ and $BD=A$. The relation 
$\mathcal L\cap \mathcal R$ is denoted by $\mathcal H$, and the 
join $\mathcal L \lor \mathcal R$ is denoted by $\mathcal D$. 
It is known that $\mathcal{L}$ and $\mathcal{R}$ commute in
every semigroup, thus we have 
$\mathcal L \lor \mathcal R = \mathcal{L} \circ \mathcal{R}$.

According to Green's theorem, the maximal subgroups of 
$\mat{n}{L}$ are precisely the $\mathcal{H}$-classes $H_E$  
of idempotent matrices $E\in \mat{n}{L}$. Moreover, 
if two idempotents $E,F$ belong to the same 
$\mathcal{D}$-class, then the groups $H_E$ and $H_F$ are
isomorphic. For further background on semigroup theory, and in
particular on Green's relations, see \cite{Howie}.

Green's relations in $\boole{n}$ can be described in terms
of in- and out-neighborhoods in the relations 
corresponding to matrices over $\two$.
We introduce the following notation for any relation
$\alpha\subseteq X^2$:
\begin{itemize}
\item $\oa(x)=\{z \mid (x,z)\in\alpha \}\subseteq X$ 
is the out-neighborhood of  $x \in X$,
\item $\oa(Y)=\{z \mid (y,z)\in\alpha 
\text{ for some } y\in Y\} = \bigcup_{y\in Y} \oa(y)$ 
is the out-neighborhood of a set $Y \subseteq X$, and
\item $\oa=\{\oa(Y) \mid Y\subseteq X\}$ 
is the set of all outneighborhoods.
\end{itemize}
The in-neighborhoods $\ia(x)$ and 
$\ia(Y)$ of vertices and of sets of vertices
and the set $\ia$ of all 
in-neighborhoods are defined dually.

\begin{remark}\label{rem:alpha+trans}
If $\alpha$ is transitive, then every element of $\oa$ is a ``forward-closed" set: $U \in \oa, u \in U, (u,x) \in \alpha$ implies $x \in U$. 
In fact, it is not hard to see that this condition is equivalent to transitivity.
If $\alpha$ is a partial order, then it may be more convenient to denote it by $\leq$, and then forward-closed sets are called \emph{upsets}.
\end{remark}

Note that $\oa$ and $\ia$ form lattices
under inclusion. The bottom element of both lattices
is $\oa(\emptyset)=\ia(\emptyset)=\emptyset$,
but the top elements of the two lattices might be different.
The join operation in  $\oa$ and in  $\ia$
is just the union, but the meet operation need not
be the intersection.

The following description of Green's relations on $\boole{n}$ is a combination of results obtained by Zaretskii in \cite{Zaretskii1963} (see also \cite{PlemWest1970}).

\begin{proposition}\label{prop:Green via alpha+}
Let $A,B\in\boole{n}$  and let
$\alpha, \beta\subseteq X^2$ be the corresponding binary 
relations. Then the following hold:
\begin{enumerate}
\item $(A,B)\in \mathcal L$ if and only if 
$\oa=\ob$;
\item $(A,B)\in \mathcal R$ if and only if 
$\ia=\ib$;
\item $(A,B)\in \mathcal H$ if and only if 
$\oa=\ob$ and $\ia=\ib$;
\item $(A,B)\in \mathcal D$ if and only if $\oa$ and 
$\ob$ are lattice isomorphic. 
\end{enumerate}
\end{proposition}

Every element of $\oa$ is a union of some sets
of the form $\oa(x)\ (x \in X)$, and likewise for
$\ia$, therefore we can reformulate
the first two items of the above proposition as follows.

\begin{proposition}\label{prop:green}
Let $A,B\in\boole{n}$  and let
$\alpha, \beta\subseteq X^2$ be the corresponding binary 
relations. Then the following hold:
\begin{enumerate}
\item $(A,B)\in \mathcal L$ if and only if 
for every $x\in X$ there exist $Y,Z\subseteq X$ such that 
$\oa(x)=\ob(Y)$ and 
$\ob(x)=\oa(Z)$ ;
\item $(A,B)\in \mathcal R$ if and only if 
for every $x\in X$ there exist $Y,Z\subseteq X$ such that 
$\ia(x)=\ib(Y)$ and 
$\ib(x)=\ia(Z)$.
\end{enumerate}
\end{proposition}

Since maximal subgroups contained in the same $\mathcal{D}$-class are 
isomorphic, we would like to find a ``simplest" idempotent in any given 
$\mathcal{D}$-class, and describe the maximal subgroup around this 
idempotent. For this, we need a way to tell whether two idempotents
are $\mathcal{D}$-related or not. The following lemma serves this purpose.

\begin{lemma}\label{lemma:DAD}
Let $A \in \boole{n}$ be an idempotent matrix and let $\alpha \subseteq X^2$ be the corresponding pseudo-order relation. 
Let $T$ be a complete system of representatives of the blocks of the equivalence relation $\alpha\cap\alpha^{-1} \subseteq \park{\alpha}^2$.
Let $D=(d_{ij})\in\boole{n}$ be the diagonal matrix defined by
\[
d_{ij}=
\begin{cases}
a_{ij},&\text{if } i=j\in T;\\
0,&\text{otherwise}.
\end{cases}
\]
Then we have $ADA=A$, and consequently $DAD$ is an idempotent matrix in the $\mathcal{D}$-class of $A$.
The pseudo-order relation corresponding to the matrix $A_1:=DAD$ is $\alpha_1:= \alpha\cap T^2$, and $\alpha_1$ is a partial order on $T$.
\end{lemma}
\begin{proof}
The entries of the matrix $ADA=(c_{ij})$ can be computed as follows (taking into account that $d_{k\ell}=0$ whenever $k\neq\ell$); for comparison we also write out the entry $a_{ij}$ from the product $A=AA$: 
\begin{align*}
c_{ij} &= \sum_{k,\ell=1}^{n}a_{ik}d_{k\ell}a_{\ell j}
= \sum_{k=1}^{n}a_{ik}d_{kk}a_{kj}
= \sum_{k=1}^{n}a_{ik}a_{kk}a_{kj}; \\
a_{ij} &= \sum_{k=1}^{n}a_{ik}a_{kj}. 
\end{align*}
It is clear that $c_{ij}\leq a_{ij}$, as $a_{ik}a_{kk}a_{kj}\leq a_{ik}a_{kj}$ for all $i,j,k\in X$. 
The inequality $c_{ij}\geq a_{ij}$ is equivalent to the implication $a_{ij}=1\implies c_{ij}=1$, and we will prove this using of the transportation network interpretation of matrices (see Remark~\ref{rem:trucks}). 
If $a_{ij}=1$ (i.e., $(i,j)\in\alpha$), then there is a (one-step) route from site $i$ to site $j$. 
Therefore, by condition (b') of Remark~\ref{rem:pseudo truck}, there is a route from $i$ to $j$ that passes through some site $p\in\park{\alpha}$ with a parking lot. 
Since $T$ is a complete system of representatives of the blocks of the equivalence relation $\alpha\cap\alpha^{-1}$, there is a (unique) $k\in T$ with $(p,k),(k,p)\in\alpha$.
By the transitivity of $\alpha$, this implies that $(i,k),(k,k),(k,j)\in\alpha$, hence $a_{ik}=a_{kk}=a_{kj}=1$. This proves that $c_{ij}=1$, thus $ADA=A$, as claimed.

Idempotence of $DAD$ now follows easily (note that $D^2=D$):
\[
(DAD)^2=DA(DD)AD=DADAD=D(ADA)D=DAD.
\]
To prove the relation $(A,DAD)\in\mathcal{D}$, we verify that $(A,DA)\in\mathcal{L}$ and $(DA,DAD)\in\mathcal{R}$:
\begin{align*}
A(DA)=A &\implies (A,DA)\in\mathcal{L};\\
(DAD)A=DA &\implies (DA,DAD)\in\mathcal{R}.
\end{align*}

If $A_1=DAD=(b_{ij})$, then $b_{ij}=d_{ii}d_{jj}a_{ij}$, thus $b_{ij}=a_{ij}$ if $i,j\in T$, and $b_{ij}=0$ otherwise. 
This means that the relation $\alpha_1$ is obtained from $\alpha$ by deleting all edges going into or out from vertices outside $T$.
The choice of the set $T$ ensures that $\alpha_1$ is a partial order on $T$ (and each element of $X\setminus T$ is an isolated vertex in $\alpha_1$). 
\end{proof}

The pseudo-order $\alpha_1$ constructed in Lemma~\ref{lemma:DAD} is a partial order on the set $T\subseteq X$. Relations of this form (i.e., partial orders on subsets of $X$) are called \emph{reduced idempotents}. It was already proved in \cite{MontPlem1969} that if a $\mathcal{D}$-class contains an idempotent, then it also contains a reduced idempotent. We complement this result with a simple criterion to decide whether two pseudo-orders are $\mathcal{D}$-related (see Theorem~\ref{thm:pseudo-D} below).
The structure of the poset $(T;\alpha_1)$ is independent
of the choice of $T$; let us denote 
(the isomorphism type of) this poset by $\poa$.
If we use the usual symbol $\leq$
for this partial order instead of $\alpha_1$, then the out- and 
in-neighborhood of $x\in T$ can be written as:
\begin{itemize}
\item $\outn{\alpha_1}(x)=\{y\in T:y\geq x\}=:\up{x}$;
\item $\inn{\alpha_1}(x)=\{y\in T:y\leq x\}=:\down{x}$.
\end{itemize}
The elements of $\outn{\alpha_1}$ (i.e., unions of sets of the form $\up{x}$) are called \emph{upsets} (cf. Remark~\ref{rem:alpha+trans}). 
Thus $U \subseteq T$ is an upset if and only if $x\in U$ and $y\geq x$ implies $y\in U$ for all $x,y\in T$.
Dually, the members of of $\inn{\alpha_1}$ are called \emph{downsets}.

\begin{theorem}\label{thm:pseudo-D}
Let $A,B \in \boole{n}$ be idempotent matrices 
and let $\alpha,\beta \subseteq X^2$
be the corresponding pseudo-order relations. 
We have $(A,B) \in \mathcal{D}$ if and only if the posets
$\poa$ and $\po{\beta}$ are isomorphic.
\end{theorem}
\begin{proof}
Let $\alpha$ be a pseudo-order on $X$, and let the set $T$ be defined as in Lemma~\ref{lemma:DAD}.
It follows from Lemma~\ref{lemma:DAD} and Proposition~\ref{prop:Green via alpha+} that the lattices $\oa$ and $\outn{\alpha_1}$ are isomorphic, and it is clear that the latter is isomorphic to the lattice of upsets of $\poa$. Thus, by Proposition~\ref{prop:Green via alpha+}, we only need to prove that the posets $\poa$ and $\po{\beta}$ are isomorphic if and only if their upset lattices are isomorphic.
The ``only if" part is trivial, and the ``if" part follows from the observation that for any poset $P$, the join-irreducible elements of the upset lattice are exactly the sets of the form $\up{x}\ (x\in P)$, and these sets form a poset that is isomorphic to $P$.
\end{proof}

\subsection{Maximal subgroups}
We conclude this section with the promised description of the maximal subgroups of $\boole{n}$. 
First we need a simple auxiliary observation.

\begin{lemma}\label{lemma:extensive}
If $(T;\leq)$ is a finite poset and $f$ is a permutation of $T$ such that $f(x)\geq x$ for all $x\in T$, then $f=\operatorname{id}_T$.
\end{lemma}
\begin{proof}
If $x$ is a maximal element, then $f(x)=x$ follows immediately from the assumption $f(x)\geq x$.
From here, we can proceed downwards, proving by induction on the size of $\up{x}$ that $f(x)=x$ for all $x\in T$.
\end{proof}

\begin{theorem}\label{thm:H-class}
Let $A\in\boole{n}$ be an idempotent matrix with the corresponding pseudo-order $\alpha\subseteq X^2$. Assume that $A$ is a reduced idempotent, i.e., $\alpha$ is a partial order on the set $T:=\park{\alpha}\subseteq X$. 
Then a matrix $B\in\boole{n}$ belongs to the $\mathcal H$-class of $A$ if and only if it can be written as $B=P_{f}A$, where $f$ is a permutation on $X$ such that $f(T)=T$ and the restriction of $f$ to $T$ is an automorphism of the poset $(T;\alpha)$ .
\end{theorem}

\begin{proof}
Let us first interpret the requirements imposed on $f$ in the theorem. 
Since all elements of $X\setminus T$ are isolated in $\alpha$, a permutation $f$ of $X$ is an automorphism of the pseudo-ordered set $(X;\alpha)$ if and only if $f(T)=T$ and $f$ is an automorphism of the poset $(T;\alpha)$. Regarding $f$ as a binary relation, it is clear that $f$ is an automorphism of $(X;\alpha)$ if and only if $f\circ\alpha\circ f^{-1}=\alpha$, which is in turn equivalent to $f\circ\alpha=\alpha\circ f$. The latter condition can be formulated in terms of matrices as $P_fA=AP_f$. If this holds, then it is easy to verify that $B=P_fA$ is $\mathcal{H}$-related to $A$:
\begin{align*}
B=P_fA \text{ and } P^{-1}_fB=A &\implies (A,B)\in\mathcal{L};\\
B=AP_f \text{ and } BP^{-1}_f=A &\implies (A,B)\in\mathcal{R}.
\end{align*}

Assume now that $B\in\boole{n}$ belongs to the $\mathcal H$-class of $A$, and let $\beta\subseteq X^2$ denote the relation corresponding to $B$. Then we have $\oa=\ob$ and $\ia=\ib$ by Proposition~\ref{prop:Green via alpha+}. If $x\in X\setminus T$, then $x$ does not appear in any member of $\oa$ or $\ia$, thus $x$ is an isolated point in $\beta$, too. Therefore, it suffices to focus on elements of $T$. 

If $x\in T$, then $x\in\oa(x)$, and this implies that if we write $\oa(x)=\ob(Y)=\bigcup_{y\in Y}\ob(y)$ for a suitable $Y\subseteq X$ as in Proposition~\ref{prop:green}, then $x\in\ob(y)$ for some $y\in Y$. Using Proposition~\ref{prop:green} again, we get a set $Z\subseteq X$ such that $\ob(y)=\oa(Z)=\bigcup_{z\in Z}\oa(z)$. Therefore, there exists $z\in Z$ with $x\in\oa(z)$, and then $\oa(x)\subseteq\oa(z)$, since $\alpha$ is transitive. Now we can conclude that $\oa(x)=\ob(y)$:
\[
\oa(x)\subseteq\oa(z)\subseteq\ob(y)\subseteq\oa(x).
\]
Note that $y\in T$, as $\oa(x)=\ob(y)$ is not empty. Letting $h(x)=y$, we can define a map $h\colon T \to T$ such that $\oa(x)=\ob(h(x))$ for all $x\in T$. If $h(x_1)=h(x_2)$, then $\oa(x_1)=\ob(h(x_1))=\ob(h(x_2))=\oa(x_2)$, and this can hold only if $x_1=x_2$, as $\alpha$ is antisymmetric. This proves that $h$ is injective, thus it is also bijective by the finiteness of $T$. We can extend the inverse of $h$ to a permutation $f$ on $X$ by keeping the elements of $X\setminus T$ fixed. The defining property $\oa(x)=\ob(h(x))$ of $h$ can be expressed in terms of $f$ as follows:
\begin{equation}\label{eq:f}
\forall x,y\in X\colon (f(x),y)\in\alpha\iff(x,y)\in\beta.
\end{equation}
(Note that the above equivalence holds trivially if $x$ or $y$ lies outside of $T$, since each element of $X\setminus T$ is an isolated vertex in $\alpha$ as well as in $\beta$.)

Repeating the previous argument for the in-neighborhoods, we obtain a permutation $g$ on $X$ such that 
\begin{equation}\label{eq:g}
\forall x,y\in X\colon (x,g(y))\in\alpha\iff(x,y)\in\beta.
\end{equation}
Comparing \eqref{eq:f} and \eqref{eq:g}, we see that
\begin{equation}\label{eq:fg}
\forall x,y\in X\colon (f(x),y)\in\alpha\iff(x,g(y))\in\alpha.
\end{equation}
In particular, since for every $x\in T$ we have $(f(x),f(x))\in\alpha$, it follows from \eqref{eq:fg} that $(x,g(f(x)))\in\alpha$ holds for all $x\in T$. 
Applying Lemma~\ref{lemma:extensive} to the restriction of the map $fg$ to $T$, we can draw the conclusion $g=f^{-1}$ (recall that both $f$ and $g$ act identically on $X\setminus T$). 
Now it is easy to deduce from \eqref{eq:fg} that $f$ is an automorphism of $\alpha$:
\[
(x,y)\in\alpha\iff(x,g(f(y))\in\alpha\iff(f(x),f(y))\in\alpha.
\]

Finally, let us observe that \eqref{eq:f} means that the matrix $B$ is obtained from $A$ by permuting its rows by the permutation $f$, hence $B=P_fA$ (see Remark~\ref{rem:permutation matrices}).
%
\end{proof}

\begin{corollary}
Let $A\in\boole{n}$ be an idempotent matrix with the corresponding pseudo-order $\alpha\subseteq X^2$. Then the $\mathcal H$-class containing $A$ is isomorphic to the automorphism group of the poset $\poa$
\end{corollary}

\begin{proof}
By Green's theorem, the $\mathcal H$-classes within the same $\mathcal{D}$-class form isomorphic groups, thus Theorem~\ref{thm:pseudo-D} (or Lemma~\ref{lemma:DAD}) allows us to assume without loss of generality that $A$ is a reduced idempotent.
According to Theorem~\ref{thm:H-class}, the elements of the $\mathcal{H}$-class $H_A$ of $A$ are of the form $P_fA$, where $f$ is a permutation on $X$ such that $f(T)=T$ and the restriction of $f$ to $T$ is an automorphism of the poset $(T;\alpha)$. 
We have seen in the proof of Theorem~\ref{thm:H-class} that $P_fA=AP_f$ holds for such permutations $f$. 
Note that the rows of $A$ indexed by elements of $X\setminus T$ are all constant zero, and the rows indexed by elements of $T$ are pairwise distinct, since $\alpha$ is antisymmetric.
Thus $P_fA$ depends only on the action of $f$ on $T$ (see Remark~\ref{rem:permutation matrices}).
Therefore, the map $\varphi:\operatorname{Aut}(\poa)\to H_A, f\mapsto P_fA$ is a well-defined bijection, and it is a group homomorphism:
\[
P_{f_1}A\cdot P_{f_2}A = P_{f_1}P_{f_2}AA = P_{f_1f_2}A. 
\]

\end{proof}

\section{Matrices over distributive lattices}\label{sec:distri}

If $L=(L;+,\cdot)$ is a bounded distributive lattice with least element $0$ and greatest element $1$, then, by Birkhoff's representation theorem, $L$ can be embedded into the lattice $\pow(\set)$ of subsets of a set $\set$ in such a way that $0$ is mapped to $\emptyset$ and $1$ is mapped to $\set$.
Identifying $L$ with its embedded image, we can actually assume that $L$ is a sublattice of $\pow(\set)$ with $0=\emptyset$ and $1=\set$.
This allows us to define a homomorphism $\cut$ from $L$ to $\two=\{0,1\}$ for each $\elem\in\set$ by
\[
\cut(a)=
\begin{cases}
1, & \text{if } \elem\in a;\\
0, & \text{if } \elem\notin a.
\end{cases}
\]
We call $\cut(a)$ the \emph{cut} of the element $a$ (at $\elem$) \cite{ZhaoJunRen2008} (also called section or zero pattern \cite{KimRoush1980}).
Since $a\subseteq\set$ is exactly the set of those elements $\elem\in\set$ for which $\cut(a)=1$, every element of $L$ is uniquely determined by its cuts.
Extending $\cut$ to matrices entrywise, we get cut homomorphisms $\cut\colon\mat{n}{L}\to\boole{n}$ for all $\elem\in\set$, and matrices are also uniquely determined by their cuts:
\begin{equation}\label{eq:cuts}
\forall A,B \in\mat{n}{L}\colon A=B \iff 
[\, \forall\elem\in\set\colon \cut(A)=\cut(B) \,].
\end{equation}

\begin{remark}\label{rem:trucks general}
Let us give an interpretation of matrices over $L$ in the spirit of Remark~\ref{rem:trucks}.
As before, we regard the elements of $X=\{ 1,\dots,n\}$ as sites (cities, store-houses, etc.), numbered from $1$ to $n$, and we think of the elements of $\set$ as different types of vehicles that can travel between these sites. 
The entry $a_{ij}\subseteq\set$ of the matrix $A\in\mat{n}{L}$ determines which vehicles can (or are allowed) to pass through the road from $i$ to $j$ (the diagonal entry $a_{ii}$ is the set of vehicles that can park at site $i$).
In other words, we have a complete directed graph on $n$ vertices, and each edge $\left(  i,j\right)  $ has a ``capacity" $a_{ij}\subseteq\set$.
(In reality, the graph is rarely complete; we can take non-existing connections into account by assigning capacity $0$.) 
Given a route $i=v_{0}\rightarrow v_{1}\rightarrow\dots\rightarrow v_{\ell}=j$ of length $\ell$, the set of vehicles that can travel all the way along this route from $i$ to $j$ is the intersection (product) of the capacities of the edges involved in the route, i.e., $a_{iv_{1}}\cdot\ldots\cdot a_{v_{\ell-1}j}$. 
We will call this element of $L$ the capacity of the route. 
The set of vehicles that can go from $i$ to $j$ on some route of length $\ell$ can be computed as the join (sum) of the capacities of the routes of length $\ell$ from $i$ to $j$, which is nothing else but the $(i,j)$-entry of $A^{\ell}$.
\end{remark}

\subsection{Idempotent matrices}
From \eqref{eq:cuts} and from the fact that each $\cut$ is a homomorphism, it follows that a matrix is idempotent if and only if all of its cuts are idempotent:
\begin{align*}
A=AA 
&\iff \forall\elem\in\set\colon\cut(A)=\cut(AA)\\
&\iff \forall\elem\in\set\colon\cut(A)=\cut(A)\cut(A).
\end{align*}
Combining this observation with Theorem~\ref{thm:idemp iff pseudo}, we get the following description of idempotent matrices over distributive lattices.

\begin{proposition}\label{prop:idemp via cuts}
A matrix $A\in\mat{n}{L}$ over a distributive lattice $L\leq\pow(\set)$ is idempotent if and only the binary relation $\alpha_\elem\subseteq X^2$ corresponding to the cut matrix $\cut(A)$ is a pseudo-order for each $\elem\in\set$.
\end{proposition}

Although Proposition~\ref{prop:idemp via cuts} certainly characterizes idempotent matrices, this characterization does not give a complete picture about the idempotent elements of the semigroup $\mat{n}{L}$, since it does not tell us which systems of pseudo-orders $\alpha_\elem\,(\elem\in\set)$ can arise as cuts of idempotent matrices.
In full generality perhaps one cannot expect a feasible solution for this problem, but for chains we can give a simple criterion.
We represent the $m$-element chain in the power set of $\set=\{ 1,\dots,m-1 \}$ as
\begin{equation}\label{eq:chain of sets}
\emptyset \subset \{ 1 \} \subset \{ 1,2 \} \subset\dots\subset \{ 1,2,\dots,m-1 \},
\end{equation}
so that the cut homomorphisms are $\cutt{1},\dots,\cutt{m-1}$.

\begin{theorem}\label{thm:idemp chain}
If $L$ is the $m$-element chain, then a matrix $A\in\mat{n}{L}$ is idempotent if and only the binary relations corresponding to the cut matrices $\cutk(A)\,(k=1,\dots,m-1)$ form a system of nested pseudo-orders $\alpha_1\supseteq\dots\supseteq\alpha_{m-1}$.
\end{theorem}

\begin{proof}
Since we represent $L$ by the chain of sets $\eqref{eq:chain of sets}$, we have the implication $k\in a\implies k-1\in a$ for all $a\in L$ and $k\in\{2,\dots,m-1\}$. This implies the inequalities $\cutt{1}(A)\geq\dots\geq\cutt{m-1}(A)$ for every matrix $A\in\mat{n}{L}$ (idempotent or not), and these inequalities translate to the containments $\alpha_1\supseteq\dots\supseteq\alpha_{m-1}$ of the corresponding relations. This together with Proposition~\ref{prop:idemp via cuts} proves the necessity of the condition formulated in the proposition.

For sufficiency, assume that we have a nested sequence of pseudo-orders $\alpha_1\supseteq\dots\supseteq\alpha_{m-1}$ on $X$. 
Define the matrix $A=(a_{ij})_{i,j=1}^{n}\in\mat{n}{L}$ by 
\[
a_{ij}=\bigl\{ k\in\{1,\dots,m-1\} : (i,j)\in\alpha_k \bigr\}.
\]
Observe that the assumed containments of the relations $\alpha_k$ guarantee that $a_{ij}$ is an element of $L$. Thus $A$ is indeed a matrix over $L$, and the binary relations corresponding to the cuts of $A$ are exactly the relations $\alpha_1,\dots,\alpha_{m-1}$. Since these are all pseudo-orders, each cut of $A$ is idempotent by Theorem~\ref{thm:idemp iff pseudo}, and then idempotence of $A$ follows from Proposition~\ref{prop:idemp via cuts}.
\end{proof}

\subsection{Green's relations and maximal subgroups}
We have proved in the previous subsection that a matrix is idempotent if and only if all of its cuts are idempotent. 
The reason behind this observation is that the definition of idempotence is simply an equality; it does not ask for the existence of certain elements. 
For ``existentially quantified" notions the situation is more complicated. 
As an example, let us consider the definition of the $\mathcal{R}$-relation:
\begin{equation*}
(A,B)\in\mathcal{R}\iff\exists C,D\in\mat{n}{L}\colon AC=B \text{ and } BD=A.     
\end{equation*}
Using the fact that the cut maps are homomorphisms, it follows that $\cut(A)\cut(C)=\cut(B)$ and $\cut(B)\cut(D)=\cut(A)$, hence $\cut(A)$ and $\cut(B)$ are $\mathcal{R}$-related in the semigroup $\boole{n}$ for all $\elem\in\set$.
However, the converse is not necessarily true. 
Given matrices $C_\elem,D_\elem\in\boole{n}$ such that $\cut(A)C_\elem=\cut(B)$ and $\cut(B)D_\elem=\cut(A)$ for all $\elem\in\set$, it is not guaranteed that there exist matrices $C,D\in\mat{n}{L}$ whose cuts are $C_\elem$ and $D_\elem$, respectively. 
In fact, an example of $\mathcal{R}$-inequivalent matrices $A,B\in\mat{2}{L}$ over the three-element chain were presented in \cite{ZhaoJunRen2008} such that both of their cuts are $\mathcal{R}$-related.

Nevertheless, as illustrated by the following theorem, in some special cases we can recover information about matrices over $L$ from their cuts.
\begin{theorem}\label{thm:H-class chain}
Let $L$ be the $m$-element chain, and let $A\in\mat{n}{L}$ be an idempotent matrix such that the binary relations $\alpha_k\subseteq X^2$ corresponding to the cut matrices $\cutk(A)\,(k=1,\dots,m-1)$ are all partial orders.
Then a matrix $B\in\mat{n}{L}$ belongs to the $\mathcal H$-class of $A$ if and only if it can be written as $B=P_{f}A$, where $f$ is a permutation on $X$ that is a common automorphism of the posets $(X;\alpha_k)\,(k=1,\dots,m-1)$.
\end{theorem}

\begin{proof}
Just as in the proof of Theorem~\ref{thm:H-class}, we can see that $f$ is an automorphism of the poset $(X;\alpha_k)$ if and only if $\cutk(P_f)$ commutes with $\cutk(A)$. 
According to \eqref{eq:cuts}, this holds for every $k$ if and only if $P_fA=AP_f$, and the latter implies that the matrix $B=P_fA$ belongs to the $\mathcal{H}$-class of $A$.

Conversely, assume that $B\in\mat{n}{L}$ is $\mathcal{H}$-related to $A$. 
Since each $\cutk$ is a homomorphism,  $(\cutk(A),\cutk(B))\in\mathcal{H}$ holds in the semigroup $\boole{n}$ for all $k\in\{1,\dots,m-1\}$.
By Theorem~\ref{thm:H-class}, for each $k$ there exists an automorphism $f_k$ of the poset $(X;\alpha_k)$ such that $\cutk(B)=P_{f_k}\cutk(A)$.
We are going to prove that $f_1=\dots=f_{m-1}$.

Let us write out \eqref{eq:f} for each cut (here $\beta_k$ denotes the binary relation corresponding to the matrix $\cutk(B)\in\boole{n}$):
\begin{equation}\label{eq:f cuts}
\forall x,y\in X\colon (f_k(x),y)\in\alpha_k\iff(x,y)\in\beta_k\quad(k=1,\dots,m-1).
\end{equation}
Since $L$ is a chain, the relations $\beta_k$ form a nested sequence (cf. the beginning of the proof of Theorem~\ref{thm:idemp chain}):
\begin{equation}\label{eq:beta nested}
\beta_1\supseteq\dots\supseteq\beta_{m-1}.
\end{equation}
For every $k\in\{ 1,\dots,m-1 \}$ and $x\in X$, we have $(f_k(x),f_k(x))\in\alpha_k$, as $\alpha_k$ was assumed to be a partial order.
Using \eqref{eq:f cuts}, this implies that $(x,f_k(x))\in\beta_k$, and then $(x,f_k(x))\in\beta_1$, by \eqref{eq:beta nested}. 
Applying \eqref{eq:f cuts} with $k=1$, we can conclude that $(f_1(x),f_k(x))\in\alpha_1$.
We can rewrite this as $(y,f_k(f_1^{-1}(y)))\in\alpha_1$ with the notation $y=f_1(x)$.
This holds for every $y\in X$, therefore $f_1=f_k$ follows from Lemma~\ref{lemma:extensive}. 

We have proved that $f:=f_1=\dots=f_{m-1}$ is a common automorphism of the posets $(X;\alpha_k)\,(k=1,\dots,m-1)$.
It remains to prove that $B=P_fA$. By \eqref{eq:cuts}, it suffices to show that the cuts of $B$ and $P_fA$ coincide:
\begin{equation*}
\cutk(B) = P_{f_k}\cutk(A) =  P_{f}\cutk(A) =  \cutk(P_f)\cutk(A) = \cutk(P_fA).   
\end{equation*}
(We used the fact that cuts preserve $0$ and $1$, hence each cut of the permutation matrix $P_f$ is itself.)
\end{proof}

\begin{corollary}
Let $L$ be the $m$-element chain, and let $A\in\mat{n}{L}$ be an idempotent matrix such that the binary relations $\alpha_k\subseteq X^2$ corresponding to the cut matrices $\cutk(A)\,(k=1,\dots,m-1)$ are all partial orders.
Then the $\mathcal H$-class containing $A$ is isomorphic to the group of common automorphisms of the posets $(X;\alpha_k)\,(k=1,\dots,m-1)$.
\end{corollary}

\begin{remark}
For chains, Theorem~\ref{thm:01irr invertible} is a special case of Theorem~\ref{thm:H-class chain}.
Indeed, if $A=I$, then each $\alpha_k$ is the equality relation on $X$, hence the group of automorphisms is $S_n$.
\end{remark}

\subsection{Nilpotent matrices}
First we give a simple criterion for the nilpotency of a matrix in terms of the underlying directed graph, and then we use it to explicitly describe nilpotent matrices over bounded distributive lattices with a meet-irreducible bottom element.
\begin{lemma}\label{lemma:nilpotent iff cycles 0}
A matrix $A\in\mat{n}{L}$ over a bounded distributive lattice $L$ is nilpotent if and only if every cycle in the directed graph corresponding to $A$ has capacity $0$.
\end{lemma}

\begin{proof}
We are going to use the interpretation of matrices outlined in Remark~\ref{rem:trucks general}.
Suppose first that every cycle has capacity $0$. 
Every route of length $n$ contains a cycle (possibly of length $1$), therefore it has capacity $0$. 
Thus we have $A^{n}=0$.

Now suppose that there is a cycle $C$ of length $\ell$ with capacity $c\neq0$.
If $\elem\in c$, then trucks of type $\elem$ can drive along $C$. Driving along $C$ several times, we see that the $(i,i)$-entry of $A^{k\ell}$ contains $\elem$ for every vertex $i$  in $C$ and for every natural number $k$.
This shows that none of the matrices $A,A^\ell,A^{2\ell},\dots$ are zero, hence $A$ cannot be nilpotent.
\end{proof}

By a strictly upper triangular matrix we mean a matrix $A\in\mat{n}{L}$ that has zeros below its main diagonal as well as on the main diagonal, i.e.,
$a_{ij}\neq0\implies i<j$.

\begin{theorem} \label{thm:0irr nilpotent}
Let $L$ be a bounded distributive lattice in which $0$ (the bottom element) is meet-irreducible.
Then a matrix $A\in\mat{n}{L}$ is nilpotent if and only if it is conjugate to a strictly upper triangular matrix, i.e., there exists a strictly upper triangular matrix $U$ and an invertible matrix $C$ such that $A=C^{-1}UC$.
\end{theorem}

\begin{proof}
If $U$ is a strictly upper triangular matrix, then we have $U\leq V$, where $V$ is the matrix having ones above the diagonal and zeros on and below the diagonal:
\begin{equation}\label{eq:matrix V}
V=
\begin{pmatrix}
0 & 1 & \dots & 1 & 1\\
0 & 0 & \dots & 1 & 1\\
\vdots & \vdots & \ddots & \vdots & \vdots\\
0 & 0 & \dots & 0 & 1\\
0 & 0 & \dots & 0 & 0
\end{pmatrix}.  
\end{equation}
In the directed graph corresponding to $V$, we have an edge from $i$ to $j$ if and only if $i<j$. 
This means that it is impossible to make a route of length $n$, hence $V^n=0$.
Since $U\leq V$, it follows that $U^n=0$, which implies that $(C^{-1}UC)^n=0$ for every invertible matrix $C$.

Conversely, let us assume that $A\in\mat{n}{L}$ is a nilpotent matrix.
Consider the relation $\alpha\subseteq X^2$ defined by $\alpha:=\{ (i,j) : a_{ij}\neq0\}$.
(If $L$ is finite, then meet-irreducibility of $0$ implies that $0$ has a unique upper cover $\elem$. 
Then $\alpha$ is the relation corresponding to the matrix $\cut(A)$, i.e., $(i,j)\in\alpha$ iff trucks of type $\elem$ are allowed to travel on the edge from $i$ to $j$.)
By Lemma~\ref{lemma:nilpotent iff cycles 0}, every cycle in the directed graph corresponding to $A$ has zero capacity, hence at least one edge of each cycle has capacity $0$, as $0$ is meet-irreducible. 
This means that $\alpha$ contains no directed cycles.
Therefore, the reflexive transitive closure of $\alpha$ is a partial order on $X$, and this partial order can be extended to a linear order $\sqsubseteq$. 
Since $\sqsubseteq$ is an extension of $\alpha$, we have
$a_{ij}\neq0\implies i\sqsubset j$ for all $i,j\in X$. 

Let $\pi$ be the permutation of $X$ given by $\pi(1)\sqsubset\dots\sqsubset\pi(n)$, and let $C=P_\pi$.
We claim that the matrix $U:=CAC^{-1}$ is strictly upper triangular. 
By the definition of the matrix $C$, we have $u_{ij}=a_{\pi(i)\pi(j)}$, hence
\begin{equation*}
u_{ij}\neq 0 \implies a_{\pi(i)\pi(j)}\neq 0 \implies \pi(i)\sqsubset\pi(j) \implies  i<j.
\end{equation*}
(The last implication is justified by the definition of $\pi$.)
Thus $U$ is indeed strictly upper triangular, and this completes the proof, as $A=C^{-1}UC$.
\end{proof}

\begin{remark}
We have seen in the proof of Lemma~\ref{lemma:nilpotent iff cycles 0} that $A\in\mat{n}{L}$ is nilpotent if and only if $A^{n}=\mathbf{0}$. 
This cannot be sharpened: the matrix $V$ given in \eqref{eq:matrix V} is nilpotent, but $V^{n-1}\neq 0$.
\end{remark}

\begin{example}
Theorem~\ref{thm:0irr nilpotent} does not necessarily remain true without the assumption on the irreducibility of $0$. 
Consider the matrix 
$A=\begin{pmatrix}
0&a\\
b&0
\end{pmatrix}$
over the lattice $(\two\times\two)\oplus\mathbf{1}$ shown in Figure~\ref{fig:2x2+1}.
It is easy to verify that $A^2=0$, but $A$ is not a conjugate of a strictly upper tranigular matrix.
Indeed, by Theorem~\ref{thm:01irr invertible}, the only invertible matrices in $\mat{2}{L}$ are the permutation matrices, hence the only conjugates of $A$ are itself and the matrix
\[
P_{(12)}^{-1}\cdot A\cdot P_{(12)}=
\begin{pmatrix}
0&1\\
1&0
\end{pmatrix}
\begin{pmatrix}
0&a\\
b&0
\end{pmatrix}
\begin{pmatrix}
0&1\\
1&0
\end{pmatrix}
=
\begin{pmatrix}
0&b\\
a&0
\end{pmatrix},
\]
and neither of them is upper triangular.
\end{example}
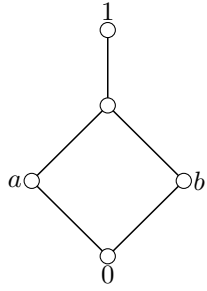
\begin{figure}
\centering
\begin{tikzpicture}[scale=1];
\draw[line width=0.2mm] (1,0) node [below]  {$0$} -- (2,1) node [right]  {$b$};
\draw[line width=0.2mm] (1,0)  -- (0,1) node [left] {$a$};
\draw[line width=0.2mm] (2,1)  -- (1,2) node [right]  {};
\draw[line width=0.2mm] (1,2)  -- (1,3) node [above]  {$1$};
\draw[line width=0.2mm] (0,1) -- (1,2);
\draw[fill=white] (1,0) circle (.1cm);
\draw[fill=white] (0,1) circle (.1cm);
\draw[fill=white] (2,1) circle (.1cm);
\draw[fill=white] (1,2) circle (.1cm);
\draw[fill=white] (1,3) circle (.1cm);
\end{tikzpicture}
\caption{The lattice $(\two\times\two)\oplus\one$}
\label{fig:2x2+1}
\end{figure}

\subsection{Fixed point iteration}
Our results on nilpotent matrices have some implications on a problem about fuzzy relations raised in \cite{SesTep2011}.
The interpretation of a matrix $A\in\mat{n}{L}$ as a directed graph with a capacity assigned to each edge (see Remark~\ref{rem:trucks general}), is almost the same as a fuzzy relation; we only need to regard the entries $a_{ij}$ as membership values instead of capacities.

Each element of a \emph{fuzzy set} has a \emph{membership value}, describing to what extent the given element belongs to the fuzzy set. 
Classically, the membership values are real numbers between $0$ and $1$, but elements of an arbitrary lattice $L$ can also serve as membership values.
(In the case of $L=\two=\{0,1\}$, we get back the usual notion of a set, called a \emph{crisp set} in this framework.) 

Thus a tuple $\x\in L^n$ describes a fuzzy subset $F$ of $X=\{1,\dots,n\}$, and a matrix $A\in\mat{n}{L}$ can be regarded as the membership function of a fuzzy subset $\alpha$ of $X^2$, which is called a \emph{fuzzy relation}.
From the definition of matrix multiplication we see that $\x A\leq\x$ is equivalent to the system of inequalities $x_i a_{ij}\leq x_j\,(i,j\in X)$, which can be interpreted as follows: ``if $i$ belongs to $F$ and $(i,j)$ belongs to $\alpha$ to some extent, then $x_j$ also belongs to $F$ at least to that extent". 
If this holds, then we say that the fuzzy set $F$ is \emph{closed} under the fuzzy relation $\alpha$.
(In the crisp case, being closed under $\alpha$ means that if $i\in F$ and $(i,j)\in\alpha$, then $j\in F$.)

The inequality $\x A\leq\x$ and the equation $\x A=\x$ were studied in \cite{SesTep2011} from the viewpoint of fuzzy control.
We refer the reader to that paper for more details about fuzzy relations and their applications, and here we focus only on the proposed fixed-point iteration method to find solutions of the equation $\x A=\x$.

The solutions of $\x A=\x$ are exactly the fixed points of the ``linear transformation" $\x\mapsto\x A$, hence we can hope that the standard fixed-point iteration method can be used to find solutions.
Thus we start with an arbitrary $\x\in L^n$, and we form the sequence 
\begin{equation}\label{eq:x,xA,xAA,...}
\x,\ \x A,\ \x A^2,\ \x A^3, \dots, \x A^k, \dots.
\end{equation}
The first problem that arises is whether this sequence converges in some topology.
Even if $L$ is an infinite lattice, each entry of each tuple in our sequence belongs to the sublattice generated by the $n+n^2$ elements $x_i,a_{ij}\,(i,j=1,\dots,n)$, which is finite if $L$ is distributive.
Therefore, the only meaningful choice is the discrete topology, and a sequence converges in this topology if and only if it becomes eventually constant.
It follows from the finiteness explained above that $\x A^k$ becomes eventually periodic.
However, the period can be longer than $1$ (consider a permutation matrix, for example), hence the sequence might fail to converge.
If $\lim_{k\to\infty}(\x A^k)$ exists, then it is easy to see that this limit will be a solution of $\x A=\x$.
Here we face a second problem: it may happen that our sequence converges to the trivial solution $\zero=(0,\dots,0)$.
To avoid this problem, it is natural start with the largest possible initial value $\one=(1,\dots,1)$.
In the following proposition we show that in this case the
limit exists, and it is the greatest solution of $\x A=\x$.

\begin{proposition}\label{prop:1A^kconverges}
Let $L$ be a bounded distributive lattice, and let $\one=(1,\dots,1)\in L^n$.
For any matrix $A\in\mat{n}{L}$, the sequence $\{\one A^k\}_{k=1}^{\infty}$ is eventually constant, and $\lim_{k\to\infty}(\one A^k)$ is the greatest solution of the fixed-pont equation $\x A=\x$.
\end{proposition}

\begin{proof}
It is clear that $\one \geq \one A$, and multiplying this inequality by $A^{k-1}$, we get $\one A^{k-1}\geq \one A^{k}$ for every natural number $k$. Therefore, we have a decreasing sequence $\one\geq\one A\geq\one A^2\geq\dots$.
Since each entry of $\one A^{k}$ belongs to the finite sublattice generated by $1$ and the entries of $A$, the sequence $\{\one A^k\}_{k=1}^{\infty}$ is eventually periodic, hence eventually constant (as it is decreasing). Thus $\one A^{\ell}=\one A^{\ell+1}=\dots$ for some natural number $\ell$, and then $\lim_{k\to\infty}(\one A^k)= A^{\ell}$.
This immediately implies that $\one A^{\ell}$ is a solution of $\x A=\x$.
Now if $\x$ is any other solution, then $\x=\x A^{\ell}\leq \one A^{\ell}$, and this means that $\one A^{\ell}$ is indeed the greatest solution.
\end{proof}

\begin{remark}
The interpetation of the greatest solution of $\x A=\x$ in the transportation network setting of Remark~\ref{rem:trucks general} is more natural if we work with column vectors instead of row vectors (or we transpose $A$). 
If $\lim_{k\to\infty}(A^k\one)=(z_1,\dots,z_n)$, then $z_i\in L$ is the set of vehicles that can start arbitrarily long trips at $i\in X$.
In other words, $z_i$ is the set of vehicles that can reach a directed cycle from $i$.
\end{remark}

Proposition~\ref{prop:1A^kconverges} allows us to completely characterize matrices $A\in\mat{n}{L}$ having a nonzero fixed point in $L^n$.

\begin{corollary}\label{cor:nonzero fixed point}
For every bounded distributive lattice $L$ and $A\in\mat{n}{L}$, the fixed-point equation $\x A=\x$ has a nonzero solution if and only if the matrix $A$ is not nilpotent.
\end{corollary}

\begin{proof}
This follows immediately from Proposition~\ref{prop:1A^kconverges}, since $\one A^{\ell}=\zero$ if and only if $A^{\ell}=0$.
\end{proof}

If the bottom element of $L$ is irreducible, then combining the above corollary with our earlier results we obtain the following corollary.

\begin{corollary}\label{cor:0irr nonzero fixed point}
Let $L$ be a bounded distributive lattice in which $0$ (the bottom element) is meet-irreducible. Then the following are equivalent for any matrix $A\in\mat{n}{L}$:
\begin{enumerate}[label=(\roman*)]
\item the only solution of the fixed-point equation $\x A=\x$ is $\zero$;
\item $\lim_{k\to\infty}(\one A^k)=\zero$;
\item $A$ is nilpotent;
\item $A^n=0$;
\item $A$ is conjugate to a strictly upper triangular matrix, i.e., there exists a strictly upper triangular matrix $U$ and an invertible matrix $C$ such that $A=C^{-1}UC$;
\item one can rearrange the rows and columns of $A$ so that it becomes a strictly upper triangular matrix, i.e., there exists a permutation $\pi\in S_n$ such that $P_\pi^{-1}AP_\pi$ is strictly upper triangular.
\end{enumerate}
\end{corollary}

Corollary~\ref{cor:0irr nonzero fixed point} applies in particular to chains (which is the most important case for applications) and it shows that the fixed-point equation $\x A=\x$ has a nontrivial solution except for a few matrices of a very restricted form.

\section{Acknowledgements}\label{sec:ack}

This research was partially supported by the National Research, Development and Innovation Office of Hungary under the FK~124814 and KH126581 funding schemes, and by grant TUDFO/47138-1/2019-ITM of the Ministry for Innovation and Technology, Hungary.

\end{document}